\documentclass[12pt, twoside, leqno]{article}
\usepackage{amsmath,amsthm}
\usepackage{amssymb,latexsym}
\usepackage{enumerate}

\newtheorem{thm}{Theorem}
\newtheorem{lemma}{Lemma}

\newtheorem*{conjcone}{Conjecture C1}
\newtheorem*{conjctwo}{Conjecture C2}
\newtheorem*{conjcthree}{Conjecture C3}
\newtheorem*{conjcfour}{Conjecture C4 {\rm (Bounded Gap Conjecture)}}

\theoremstyle{definition}

\numberwithin{equation}{section}

\allowdisplaybreaks

\def\beq{\begin{equation}}
\def\eeq{\end{equation}}
\def\cite#1{{\rm [#1]}}

\begin{document}

\title{On the difference of primes}

\author{by\\
J\'anos Pintz}

\date{}

\maketitle

\renewcommand{\thefootnote}{}

\footnote{Supported by OTKA (Grants No. K72731, K67676) and ERC/AdG. 228005.}


\renewcommand{\thefootnote}{\arabic{footnote}}
\setcounter{footnote}{0}



\section{Introduction}
\label{sec:1}
In the present work we investigate some approximations to generalizations of the twin prime conjecture. The twin prime conjecture appeared in print already the first time in a more general form, due to de Polignac \cite{Pol} in 1849:

\begin{conjcone}
Every even number can be written in infinitely many ways as the difference of two consecutive primes.
\end{conjcone}

Kronecker \cite{Kro} mentioned in 1901 the same conjecture in a weaker form as

\begin{conjctwo}
Every even number can be expressed in infinitely many ways as the difference of two primes.
\end{conjctwo}

Finally, Maillet \cite{Mai} formulated it in 1905 as

\begin{conjcthree}
Every even number is the difference of two primes.
\end{conjcthree}

We remark that it is easy to see that the existence of at least one even number satisfying either C1 or C2 is equivalent to

\begin{conjcfour}
If $p_n$ denotes the $n$\textsuperscript{th} prime then
\beq
\label{eq:1.1}
\liminf_{n \to \infty} (p_{n + 1} - p_n) < \infty.
\eeq
\end{conjcfour}

We will therefore concentrate on Conjecture~C3 which shows the closest analogy to Goldbach's conjecture among Conjectures C1--C3.

Let us call $n$ a Goldbach number if it is the sum of two primes and a Maillet number if it is the difference of two primes.
As an approximation to the conjectures of Goldbach and Maillet (C3) one can ask how long an interval can be if it contains no Goldbach numbers (i.e.\ no sum of two primes) or, no Maillet numbers (i.e.\ no difference of two primes), respectively.
Montgomery--Vaughan \cite{MV} and Ramachandra \cite{Ram} observed that if for some $0 \leq \vartheta_1, \vartheta_2 \leq 1$
\beq
\label{eq:1.2}
\pi \bigl(x + x^{\vartheta_1}\bigr) - \pi(x) \geq \frac{cx^{\vartheta_1}}{\log x}
\eeq
and
\beq
\label{eq:1.3}
\pi \bigl(n + n^{\vartheta_2}\bigr) - \pi(n) > 0 \ \text{ for almost all }\ n \in [x, 2x) \ \text{ where }\ x \to \infty
\eeq
then for $x > x_0$ the interval
\beq
\label{eq:1.4}
I = \bigl[ x, x + C'x^{\vartheta}\bigr], \quad \vartheta = \vartheta_1\vartheta_2
\eeq
contains at least one Goldbach number.
The sharpest known results, $\vartheta_1 = 21/40$ by R. C. Baker, G. Harman and the author \cite{BHP} and $\vartheta_2 = 1/20$ of Ch.\ Jia \cite{Jia} imply that $I$ contains Goldbach numbers if $\vartheta = 21/800$.
The same method applies to Maillet numbers without any change.
Thus the interval
\beq
\label{eq:1.5}
I^* = [x, x + Cx^{21/800}], \quad x > x_0
\eeq
contains an even integer which can be written as the difference of two primes.

Under supposition of the Riemann Hypothesis (RH) it was proved by Linnik \cite{Lin}, later by K\'atai \cite{Kat} that the interval
\beq
\label{eq:1.6}
I(C_2) = \bigl[ x, x + C_1(\log x)^{C_2}\bigr]
\eeq
contains Goldbach numbers for $C_2 > 3$ \cite{Lin}, respectively for $C_2 = 2$ \cite{Kat}.

We announce the unconditional

\begin{thm}
\label{th:1}
A positive proportion of even numbers in an interval of type $\bigl[x, x+ (\log x)^C\bigr]$ can be written as the difference of two primes if $C>C_0$ and  $x>x_0$.
\end{thm}

In order to illustrate the method we will prove here a result which shows that the best known exponent $21/800$ (see \eqref{eq:1.5}) can be replaced by an arbitrary positive number.

\begin{thm}
\label{th:2}
Let $\varepsilon > 0$ be arbitrary.
The interval $[x, x + x^\varepsilon]$ contains even numbers which can be written as the difference of two primes if $x > x_0(\varepsilon)$.
\end{thm}

It was shown in \cite{GPY} that the Bounded Gap Conjecture C4 is true, equivalently there is at least one de Polignac number if primes have an admissible level $\vartheta > 1/2$ of distribution.
This means that
\beq
\label{eq:1.7}
\sum_{q \leq x^{\vartheta - \varepsilon}} \max_{\substack{a\\ (a, q) = 1}} \biggl|\sum_{\substack{p \equiv a (\text{\rm mod}\, q)\\
p \leq x}} \log p - \frac{x}{\varphi(q)} \biggr| \ll_{\varepsilon, A} \frac{x}{(\log x)^A}
\eeq
for any $\varepsilon > 0$ and $A > 0$.
In \cite{Pin1} it was proved that if $\vartheta > 1/2$ then de Polignac numbers have a positive (lower) density.
We announce here the stronger but still conditional

\begin{thm}
\label{th:3}
If primes have an admissible level $\vartheta > 1/2$ of distribution then there exists a constant $C(\vartheta)$ such that for $x > C_0(\vartheta)$ the interval $[x, x + C(\vartheta)]$ contains at least one number which can be written in infinitely many ways as the difference of two consecutive primes.
\end{thm}

About 60 years ago Erd\H{o}s \cite{Erd} and Ricci \cite{Ric} independently proved that the set $J$ of limit points of the sequence
\beq
\label{eq:1.8}
\frac{p_{n + 1} - p_n}{\log p_n}, \ \text{ or equivalently that of } \ \frac{p_{n + 1} - p_n}{\log n}
\eeq
has positive Lebesgue measure but no finite point of $J$ was known until \cite{GPY}, which implied $0 \in J$.

Supposing $\vartheta > 1/2$ we can show the much stronger

\begin{thm}
\label{th:4}
Let $g(n) < \log n$ be any monotonically increasing positive function with $\lim\limits_{n \to \infty} g(n) = \infty$.
Let us suppose that primes have an admissible level $\vartheta > 1/2$ of distribution.
Then we have a constant $c(g, \vartheta)$ such that
\beq
\label{eq:1.9}
[0, c(g, \vartheta)] \subset J.
\eeq
\end{thm}

The details of proofs for Theorems~\ref{th:1}, \ref{th:3} and \ref{th:4} will appear elsewhere.
We remark here that their proofs (similar to that of Theorem~\ref{th:2}) will be completely ineffective (independently of some eliminable ineffectivity originating from the use of Bombieri--Vinogradov theorem, which uses the ineffective theorem of Siegel for $\mathcal L$-zeros).

\section{Notation and lemmata}
\label{sec:2}

Let $\mathcal P$ denote the set of primes.
Let $\varepsilon > 0$ and let $C_0$ be a sufficiently large constant depending on $\varepsilon$.
Let us suppose the existence of an infinite sequence $(H_\nu \in \mathbb Z)$
\beq
\label{eq:2.1}
I_\nu = [H_\nu, H_\nu + H_\nu^{\varepsilon}], \quad H_\nu^{\varepsilon} > 2 H_{\nu - 1}, H_1 > C_0
\eeq
such that
\beq
\label{eq:2.2}
(\mathcal P - \mathcal P) \cap \biggl(\bigcup\limits_{\nu = 1}^\infty I_\nu\biggr) = \emptyset.
\eeq

Let
\beq
\label{eq:2.3}
k = \lceil (6/\varepsilon)^2\rceil.
\eeq

We call a $k$-tuple $\mathcal H = \{h_i\}^k_{i = 1}$, $0 \leq h_1 < h_2 < \dots < h_k$ $(h_i \in \mathbb Z)$ admissible if the number $\nu_p (\mathcal H)$ of residue classes covered by $\mathcal H$ $\text{\rm mod}\, p$ satisfies
\beq
\label{eq:2.4}
\nu_p(\mathcal H) < p \ \text{ for } \ p \in \mathcal P.
\eeq
In this case the corresponding singular series is positive, i.e.\
\beq
\label{eq:2.5}
\mathfrak S(\mathcal H) := \prod_p \left(1 - \frac1p\right)^{-k} \left(1 - \frac{\nu_p(\mathcal H)}{p}\right) > 0.
\eeq
Since \eqref{eq:2.4} is trivially true for $p > k$ it is easy to choose an admissible system
\beq
\label{eq:2.6}
\mathcal H_k = \{h_\nu\}^k_{\nu = 1}, \ \ h_\nu \in I'_\nu := \bigl[ H_\nu + H_\nu^\varepsilon / 2, \ H_\nu + H_\nu^\varepsilon\bigr]
\eeq
if $C_0$ was chosen large enough (depending on $\varepsilon$).
This system $\mathcal H_k = \mathcal H$ will be fixed for the rest of the work.

We now choose any sufficiently large $t > \nu_0(\varepsilon, C_0)$ and denote
\beq
\label{eq:2.7}
T := H_t^{\varepsilon} / 2, \ \ N: = H_t + T, \ \ I'_t := [N, N + T] \subset I_t.
\eeq

It will be very important that our condition \eqref{eq:2.1} guarantees that for $h_\nu \in I_\nu$, $h_\mu \in I_\mu$, $\nu < \mu$ we have
\beq
\label{eq:2.8}
h_\mu - h_\nu \geq H_\mu + H_\mu^\varepsilon /2 - H_{\mu - 1} \geq H_\mu, \ \text{ so } \ h_\mu - h_\nu \in I_\mu .
\eeq
We will use this relation for $\nu, \mu \in [1, k] \cup \{t\}$.

We will choose our weights $a_n$ depending on $\mathcal H = \mathcal H_k$ exactly as in \cite{GPY}:
\beq
\label{eq:2.9}
a_n := \Lambda_R(n; \mathcal H, \ell)^2 := \Biggl( \frac1{(k + \ell)!} \sum_{\substack{ d \mid P_{\mathcal H}(n)\\ d \leq R}} \mu(d) \left(\log \frac{R}{d}\right)^{k + \ell}\Biggr)^2,
\eeq
where we define
\beq
\label{eq:2.10}
P_{\mathcal H}(n) := \prod^k_{i = 1} (n + h_i), \ \ell = \left[\frac{\sqrt{k}}{2}\right], \ \mathcal L = \log N, \ R = N \exp\left(-\sqrt{\log N}\right).
\eeq
Let further $\chi_{\mathcal P}(n)$ denote the characteristic function of the primes.
The notation $n \sim N$ will abbreviate $n \in [N, 2N]$.

We remark that for any admissible $k$-tuple $\mathcal H$ we have
\beq
\label{eq:2.11}
\mathfrak S(\mathcal H) \geq \prod_{p \leq 2k} \frac1p \prod_{p > 2k} \left(1 - \frac{k}{p}\right) \left(1 - \frac1p\right)^{-k} = c_1(k).
\eeq
In the following formulae the $o$ symbol will refer to the case $N \to \infty$ which is equivalent to $t \to \infty$.
We will use the notation
\beq
\label{eq:2.12}
A := \sum_{n \sim N} a_n, \ \ B := B_R(N, k, \ell) := {2\ell \choose \ell} \frac{\log^{2k + \ell} R}{(k + 2\ell)!}.
\eeq
The following lemmata are special cases of Propositions~1 and 2 of \cite{GPY}. The only change is that although Proposition~2 has originally a condition $h_i \leq R$, this can in fact be replaced without any change in the proof by $h_i \leq 4N$, for example.
Lemma~2.4 is a well-known sieve estimate (see Theorem~4.4 of \cite{Mon}, for example).

\begin{lemma}
\label{lem:2.1}
$A = \bigl(\mathfrak S(\mathcal H) + o(1)\bigr) B$.
\end{lemma}

\begin{lemma}
\label{lem:2.2}
If $h_i \in \mathcal H$ then
\beq
\label{eq:2.13}
S_i := \sum_{n \sim N} a_n \chi_{\mathcal P}(n + h_i) = \frac{2\ell + 1}{2\ell + 2} \frac{\bigl(\mathfrak S(\mathcal H) + o(1)\bigr)B}{k + 2\ell + 1}.
\eeq
\end{lemma}

\begin{lemma}
\label{lem:2.3}
If $h_0 \notin \mathcal H$, $h_0 \leq 4N$ then
\beq
\label{eq:2.14}
S_0 := S(h_0) := \sum_{n \sim N} a_n \chi_{\mathcal P} (n + h_0) = \frac{\bigl(\mathfrak S(\mathcal H\cup \{h_0\}) + o(1)\bigr)B}{\mathcal L}.
\eeq
\end{lemma}

\begin{lemma}
\label{lem:2.4}
$\pi(x + y) - \pi(x) \leq \frac{2y}{\log y}$ for any $x, y \geq 1$.
\end{lemma}

Our last lemma is Lemma~3.4 of \cite{Pin3} (and a slight variant of Theorem~2 of \cite{Pin2}).

\begin{lemma}
\label{lem:2.5}
Let $\mathcal H_k \subset [0, H]$ be a fixed $k$-element admissible set, $k \geq k_0$, $C$ a sufficiently large constant, $m \in \mathbb Z$, $\delta > 0$.
Then, for $H > \exp \left(\frac{Ck}{\delta \log k}\right)$ we have
\beq
\label{eq:2.15}
S_{\mathcal H} (M, H) = \frac1H \sum_{m \in [M, M + H]} \frac{\mathfrak S(\mathcal H \cup \{m\})}{\mathfrak S(\mathcal H)} \geq 1 - \delta.
\eeq
\end{lemma}

\section{Proof of Theorem 2}
\label{sec:3}

We begin with the very important remark that (cf.\ \eqref{eq:2.2} and \eqref{eq:2.8}) due to our definition of $I$ and to $h_\mu - h_\nu \in I_\mu$ (for $\mu > \nu$), the set $\{n + h_i\}^k_{i = 1}$ contains at most one prime for any given~$n$.
Let
\beq
\label{eq:3.1}
\mathcal D_0 = \{n \sim N; \ n + h_i \notin \mathcal P, \ i = 1,2,\dots, k\}, \ \ \mathcal D_1 = [N, 2N] \setminus \mathcal D_0.
\eeq
Then we have by \eqref{eq:2.10}--\eqref{eq:2.12} and Lemmas~\ref{lem:2.1} and \ref{lem:2.2}, for $k > k_0$
\begin{align}
\label{eq:3.2}
A_0 :&= \sum_{n \in \mathcal D_0} a_n = A - \sum^k_{i = 1} S_i \\
&= \bigl(\mathfrak S(\mathcal H) + o(1)\bigr) B \left(1 - \left(1 - \frac1{2\ell + 2}\right) \left(1 - \frac{2\ell + 1}{k + (2\ell + 1)}\right)\right)\notag\\
&\leq \frac{7 \mathfrak S(\mathcal H)B}{3\sqrt{k}}.\notag
\end{align}

The second important remark is that if $n \in \mathcal D_1$ and $m \in I'_t$ then $n + m \notin \mathcal P$.
Namely, $n \in \mathcal D_1$ and $n + m \in \mathcal P$ would imply $n + h_i \in \mathcal P$ for some $i$; consequently
\beq
\label{eq:3.3}
m - h_i \in \mathcal P - \mathcal P , 
\eeq
which is a contradiction to \eqref{eq:2.2} in view of $m \in I_t$.
Hence,
\beq
\label{eq:3.4}
S := \sum_{n \in \mathcal D_0} a_n \sum_{m \in I'_t} \chi_{\mathcal P} (n + m) = \sum_{m \in I_{t'}} \sum_{n \sim N} a_n \chi_{\mathcal P}(n + m).
\eeq

Considering the LHS, from \eqref{eq:3.2}, \eqref{eq:2.3} and Lemma~\ref{lem:2.4} we obtain
\beq
\label{eq:3.5}
S \leq A_0 \cdot \frac{2T}{\log T} < \frac{5\mathfrak S(\mathcal H) BT}{\varepsilon \mathcal L \sqrt{k}} \leq \frac{5\mathfrak S(\mathcal H)BT}{6\mathcal L}.
\eeq
On the other hand, considering the RHS, Lemmas \ref{lem:2.3} and \ref{lem:2.5} imply
\beq
\label{eq:3.6}
S \geq \sum_{m \in I_{t'}} \frac{\bigl(\mathfrak S(\mathcal H \cup \{m\}) + o(1)\bigr)B}{\mathcal L} > \frac{5\mathfrak S(\mathcal H)BT}{6\mathcal L},
\eeq
which contradicts \eqref{eq:3.5}. This proves Theorem~\ref{th:2}.

\bigskip
\noindent
J\'anos Pintz\\
R\'enyi Mathematical Institute of the Hungarian Academy
of Sciences\\
Budapest, Re\'altanoda u. 13--15\\
H-1053 Hungary\\
E-mail: pintz@renyi.hu

\end{document}